\documentclass[12pt,a4paper]{article}
\usepackage[T1]{fontenc}
\usepackage{lmodern}
\usepackage[latin1]{inputenc}
\usepackage{latexsym}                 
\usepackage{color}                 
\usepackage{amssymb}
\usepackage{amsmath}
\usepackage{amsthm}
\usepackage{graphicx}
\usepackage{hyperref}
\usepackage{verbatim}
\usepackage{mathptmx}      
\def\eref#1{(\ref{#1}%
)}
\def\RSref#1{\ref{#1}%
}
\def\RSlabel#1{\label{#1}%
}
\def\RScite#1{\cite{#1}%
}

\newcommand{\bql}[1]{%
\begin{equation}\label{#1}%
}

\def\filename#1{}

\newcommand{\eq}{\end{equation}}
\def\dfrac#1#2{\displaystyle{\frac{#1}{#2}   }}

\def\b1{\mathbf 1}

\newcommand{\R}{\ensuremath{\mathbb{R}}}

\def\biglf{\par\bigskip\noindent}
\newtheorem{definition}{Definition}

\newtheorem{theorem}{Theorem}

\usepackage{amsmath}
\usepackage{upgreek}
\def\biglf{\par\bigskip\noindent}
\begin{document}
\begin{center}
{\large \bf $H$-Sets for Kernel-Based Spaces}
\biglf
Robert Schaback\footnote{%
Institut für Numerische und Angewandte Mathematik,
Universit\"at G\"ottingen, Lotzestra\ss{}e 16--18, 37083 G\"ottingen, Germany,
{\tt schaback@math.uni-goettingen.de}}\\
\biglf
Draft of \today
\end{center}
    \biglf
    \biglf
{\bf Abstract: } The concept of $H$-sets as introduced by Collatz in 1956
was very useful in univariate Chebyshev approximation
by polynomials or Chebyshev spaces. In the multivariate setting,
the situation is much worse, because there is no alternation, and $H$-sets
exist, but are 
only rarely accessible by mathematical arguments.
However, in Reproducing Kernel Hilbert spaces,
$H$-sets are shown here to have  a rather simple and complete characterization.
As a byproduct, the strong connection of $H$-sets to Linear Programming
is studied.
But on the downside, it is explained why $H$-sets
have a very limited range of applicability in the times of large-scale
computing.
\section{$H$-Sets and Their Use}\RSlabel{SecIntro}
Let $F$ be a space of continuous real-valued
functions on a compact domain $T$,
and consider linear approximations of functions $f$ by functions $v$ from
a subspace $V$ of $F$. In 1956, Lothar
Collatz \RScite{collatz:1956-1} introduced
\begin{definition}\RSlabel{DefHsets}
An $H$-set for $V\subseteq F\subseteq C(T)$
consists of a subset $H$ of $T$ and a sign function $\sigma\;:\;H\to \{-1,+1\}$
such that there is no $v\in V$ that makes
all values $v(h)\sigma(h)$
for $h \in H$ negative.
\end{definition} 
The classical application is in linear Chebyshev approximation
\RScite{collatz:1956-1},
stated here
in abstract form:
\begin{theorem}\RSlabel{TheAlt}
Assume that a user
has found some candidate $\tilde v\in V$ for approximation of $f\in F$
by functions from $V$,
and an $H$-set
consisting of $H$ and $\sigma$. If furthermore
\bql{eqmu}
\inf_{h\in H}(f(h)-\tilde v(h))\sigma(h)=:\mu
\eq
is positive, then 
$$
\mu \leq \displaystyle{  \inf_{v\in V}\|f-v\|_\infty}\leq \|f-\tilde v\|_\infty
$$
bounds the optimal approximation error from both sides by observable quantities.
\end{theorem}
{\bf Proof}:
For any $v\in V$, the expression
$$
\begin{array}{rcl}
  \|f-v\|_\infty
  & \geq &
  (f(h)-v(h))\sigma(h) \\
  &=&
  (f(h)-\tilde v(h))\sigma(h)+ (\tilde v(h)-v(h))\sigma(h)  \\
\end{array}
$$
implies
$$
\|f-v\|_\infty\geq (f(\tilde h)-\tilde v(\tilde h))\sigma(\tilde h)\geq \mu
$$
for some $\tilde h\in H$. $\Box$
\biglf
This shows that $H$-sets should pick near-extremal points of the error
function and keep the sign of the error there. In 1956, 
computations were still made mechanically, and then $H$-sets allowed to
assess the quality of an approximation without any large-scale computation.
\biglf
If $H$ has only $N$ points, if $V$ is $n$-dimensional, and if
the corresponding discrete Chebyshev approximation on $H$ is carried out
exactly by a linear optimizer of Simplex type,
one gets an $H$-set based on extremal points for free, as we shall prove
in Theorem 
\RSref{Thebestapp} below.
However, Theorem \RSref{eqmu} is useless in that case, because
best approximation errors on $T$ always have lower bounds
by best approximation errors on subsets. This implies that
the merits of specially constructed $H$-sets are restricted to inexact
discrete Chebyshev approximation. 
\section{Examples}\RSlabel{SecExa}
The simplest and classical example is Chebyshev approximation in $C[-1,+1]$
by polynomials of degree $n$. One can expect alternation of the error
of best Chevyshev approximations on sets $T$ 
of $n+2$ points, and these are the canonical candidates for an $H$ set,
the signs being alternating wrt. the ordering of the points. The Remes exchange
algorithm makes heavy use of this principle, and Theorem
\RSref{TheAlt} can be applied as soon as sign patterns and extremal points
stabilize in the iteration. 
\biglf
For multivariate approximation, there is no alternation principle,
and $H$-sets may be very hard to determine. But we shall
see in section \RSref{SecKC} that this is not the case for kernel-based spaces.
\biglf
In general, after performing some numerical approximation,
one may choose near-extremal points, with signs related to the sign of the error
there to get a candidate for an $H$-set, satisfying \eref{eqmu}, but then 
one must hope for the 
$H$-set property for that choice of signs. 
\biglf
Conversely, one might prove the $H$-set property for a fixed choice of $V$,
$H$, and
$\sigma$, but then the application requires these signs to arise
in \eref{eqmu}, limiting the applicability seriously.
\biglf
This gap is a general obstacle
to the practical applicability of $H$-sets.
\section{Connections to Linear Optimization}\RSlabel{SecDA}
It is strange that  most of the literature on $H$-sets
(see e.g. Taylor \RScite{taylor:1972-1}, 
Brannigan \RScite{brannigan:1977-1}, Dierieck \RScite{dierieck:1977-1}, and 
Brannigan \RScite{brannigan:1983-1})
focuses on minimality and geometry of $H$-sets and has some links
to duality, but no explicit connection to Linear Optimization.
The only exception seems to be Wetterling \RScite{wetterling:1979-1}
who briefly mentions the connection of $H$-sets to the dual Simplex algorithm.
We give details here, to prepare for the kernel-based case.
\biglf
If $V$ is $n$-dimensional with basis
$v_1,\ldots,v_n$ and if $H$ has $N$ points $h_1,\ldots,h_N$
with associated signs $\sigma_1,\ldots,\sigma_N$,
one can form the $N\times n$ matrix $A$ with entries $v_i(h_k)\sigma_k$.
\begin{theorem}\RSlabel{ThedualH}
  Under the above notation, the $H$-set property is equivalent to
  the two equivalent dual statements:
  \begin{itemize}
 \item There is no $x\in\R^n$ such that the vector
   $b:=Ax\in \R^N$ is negative in all components,
 \item There is a nonzero nonnegative vector $w\in\R^N$ with
\bql{eqdual}
w^TA=0= \sum_{k=1}^Nw_kv_i(h_k)\sigma_k,\;1\leq i\leq n.
\eq
\end{itemize} 
\end{theorem}
{\bf Proof}: The first statement is the definition of the $H$-set property.
The second implies the first, because $w^TAx=0=w^Tb$ makes it impossible
that $b:=Ax\in \R^N$ is negative in all components.
The converse is also true, due to the Farkas lemma in the background:
\begin{itemize}
 \item[] $Ax\leq b$ is solvable if and only if for all vectors $w\geq 0$ and $w^TA=0$
the inequality $w^Tb\geq 0$ holds.
\end{itemize}
If we have an $H$-set, the problem $Ax\leq -\beta 1$
is unsolvable for small fixed $\beta>0$. This implies that $Ax\leq -1$ is
unsolvable, and then there is a $w\geq 0$ with $w^TA=0$ and  $-1^Tw<0$.
$\Box$
\biglf
We shall use \eref{eqdual} for a 
numerical test for the $H$-set property.
To decide that $H$ and $\sigma$ form an $H$-set or not, we pose the
solvable problem
\bql{eqPro1w}
\begin{array}{rcl}
  1^Tw&=& Max!\\
  0\leq& w&\leq 1\\
  A^Tw&=&0,
\end{array}
\eq
start at the origin and check if the maximum is positive or zero.
\biglf
The condition \eref{eqdual} means that there is a point evaluation
functional based on $H$ that vanishes on $X$, and the signs are determined
by the coefficients of the functional. This is very useful
when approximating  with univariate polynomials of degree $n$ on $n+2$ points,
because the required functional is the divided difference up to a factor.
In general, the signs have a dual role: they arise in a primal sense
as signs of function values
and in a dual sense as signs of coefficients of functionals.
The duality is twofold: {\it values $\Leftrightarrow$ coefficients} and
{\it functions $\Leftrightarrow$ functionals}.
\biglf
There is another connection to Linear Optimization that explains why
$H$-sets lost much of their importance in presence of large-scale computing.
This elaborates a short remark by
Wetterling \RScite{wetterling:1979-1}.
\begin{theorem}\RSlabel{Thebestapp}
  If best discrete Chebyshev approximation in finite-dimensional spaces
  is  written as a Linear Optimization problem, one gets an $H$-set
  as a subset of extremal points with associated signs for free,
  provided that calculations are exact and a solution of the dual
  problem is provided as well.
\end{theorem} 
{\bf Proof}: 
For discrete Chebyshev approximation of data $f_H\in\R^N$
on $H$, using the $N\times n$ matrix $B$ with entries $v_i(h_k)$,
one can pose the
linear optimization problem 
\bql{eqprimalCAP0}
\begin{array}{rcl}
  \eta&=& Min!\\
  \left (\begin{array}{rcl}
-B & -1_H\\B & -1_H
  \end{array} \right)
  \left(\begin{array}{rcl}
x\\ \eta
\end{array} \right)
  &\leq& \left(
  \begin{array}{rcl}
-f_H\\f_H
  \end{array}
  \right)
\end{array}
\eq
and the dual problem is to find some $w\in R^N$ with 
$$
\begin{array}{rcl}
  f_H^Tw&=& Max!\\[0.1cm]
  B^Tw&=& 0\\[0.1cm]
  \|w\|_1&=&1
\end{array}
$$
to be implemented via a split $w=w^+-w^-$ in positive and negative parts.
Both problems
are solvable, and if $w^*,\,x^*$, and $\eta^*$ are the optimal solutions,
one has
$$
\begin{array}{rcl}
  f_H^Tw\leq f_H^Tw^*&=&\eta^*=\|f_H-Bx^*\|_{\infty,H}\leq
  \|f_H-Bx\|_{\infty,H}\\
  B^Tw^*&=& 0\\
  \sigma^*&:=& sgn(f_H-Bx^*)\\
  w_k^*&=& 0 \;\;\hbox{ if }\;\; |f-Bx^*|_k<\eta^*\\
  sgn(w_k^*)&=& \sigma^*_k \hbox{  or }  w_k^*=0\hbox{ otherwise } \\
\end{array}
$$
due to strong duality and complementary slackness.
Therefore the support of $w^*$, being a subset of the extremal
points, forms an $H$-set for free. This assumes that the
optimizer for \eref{eqprimalCAP0} is exact and 
provides the dual solution, but modern interior point
methods may fail to do so. $\Box$
\biglf
For $N$ much larger than $n$, there may be many choices of $H$-sets.
The cited literature considers minimal $H$-sets at length.
In view of minimality,
the above formulation provides $H$ sets that not necessarily have
a minimal number of points, but the minimal sum of positive weights
in the dual solution vector $w$. By use of the 1-norm, chances are good
that the optimization concentrates weights into few nonzero components,
and this can be observed in the example below. 
\section{The Kernel Case}\RSlabel{SecKC}
We now apply this to kernel-based spaces
and use the inherent duality
principles there.
Readers are referred to books
\RScite{buhmann:2004-1, wendland:2005-1,fasshauer-mccourt:2015-1}
for the background.
\biglf
Let $K$ be a symmetric
strictly positive
definite kernel on $T$, and let $V_X$ be spanned by translates
$K(\cdot,x_1),\ldots,K(\cdot,x_n)$ for $n$ different points
$x_1,\ldots,x_n$ in $T$ forming a set $X$.
The candidates for $H$-sets consist of
points $h_1,\ldots,h_N$ in $T$ forming a set $H$,
with associated signs $\sigma_1,\ldots,\sigma_N$. This also defines a
subspace $V_H$
of $F$ spanned by the $H$-translates of the kernel.
\begin{theorem}\RSlabel{TheHsets}
  The $H$-sets for $V_X$ based on a finite point set $H$ of $N$ points
  are completely
  characterized by nonzero functions $f$ in $V_H$ that vanish
  on $X$, with signs
  of the coefficients of $f$ in the basis of $V_H$.
\end{theorem} 
{\bf Proof:}
In the kernel case, the $N\times n$
matrix $A$ of the duality argument in section \RSref{SecDA}
has entries $K(x_i,h_k)\sigma_k$. Consequently,
the $H$-set property is equivalent to existence of a nonnegative  nonzero
vector $w\in\R^N$  such that
\bql{eqwKs}
\sum_{k=1}^Nw_kK(x_i,h_k)\sigma_k=0,\;1\leq j\leq n
\eq
proving the assertion for 
$$
f(x)=\sum_{k=1}^Nw_kK(x,h_k)\sigma_k.\qed
$$
Surprisingly, Theorem \RSref{TheHsets} gives a simple
characterization of all $H$-sets
in the kernel-based case, avoiding Linear Optimization completely.
Kernel spaces allow to rephrase the functional
of section \RSref{SecDA} in terms of a function.
They remove the {\it function} $\Leftrightarrow$ {\it functional} duality,
but not the {\it values} $\Leftrightarrow$ {\it coefficients} duality.
\biglf
Even in kernel-based spaces there is no nice connection
of signs of coefficients to signs of values. The only exceptions
known so far are generated by eigenvectors of kernel matrices.
There, values are positive multiples of coefficients. 
\biglf
A simple illustrative case is where $f$ is in $V_H$, and $s_{X,f}\in V_X$
interpolates $f$ on a subset 
$X$ of $H$.
Then $f-s_{X,f}$ vanishes on $X$ and determines a candidate for an $H$-set,
but then the signs of the coefficients of $f-s_{X,f}$  in the basis of the
$K(\cdot, h_k)$ should be the signs of the
values of $f-s_{X,f}$  on $H$. Such a correspondence of signs of values and
coefficients could only be expected if kernel matrices and their inverses
were sign-preserving.
\biglf
If, during a numerical approximation, the set $H$ is chosen by extremal points,
with signs determined by the error there, it is not guaranteed that
there is a function based on $H$ that vanishes on $X$ and has the required  signs of
coefficients. Conversely, if Theorem \RSref{TheHsets} can be used,
it can only be applied to cases where \eref{eqmu} has the correct signs of the
error. Even in kernel-based spaces, this gap cannot be bridged.
\section{Numerical Example}\RSlabel{SecNumExp}
Lothar Collatz always insisted that papers should have a numerical example.
Let the kernel be the Gaussian at scale one, and choose 25 points at random
in $[-1,+1]^2$ to define $X$ and the approximating space $V_X$ of translates of
the Gaussian. Then approximate the MATLAB peaks function on a regular
set $T$ of 11x11=121 points in $[-1,+1]^2$. The Chebyshev error
on $T$ comes out to be 0.0768,
while we get 0.1053 on a 41x41 evaluation grid.
The interior point method {\it lipsol}  within MATLAB's {\it linprog}
fails fo yield $H$-sets under various circumstances,
in contrast to Theorem \RSref{Thebestapp}.
If, for instance,  Lagrange multipliers larger than 1.e-5 are used,
39 points are selected with $\mu=0.0596$,
see Figure \RSref{FigLag}. Testing the $H$-set property was
done by solving the problem \eref{eqPro1w}.
\begin{figure}[hbt]
  \begin{center}
    \includegraphics[width=10.0cm,height=10.0cm]{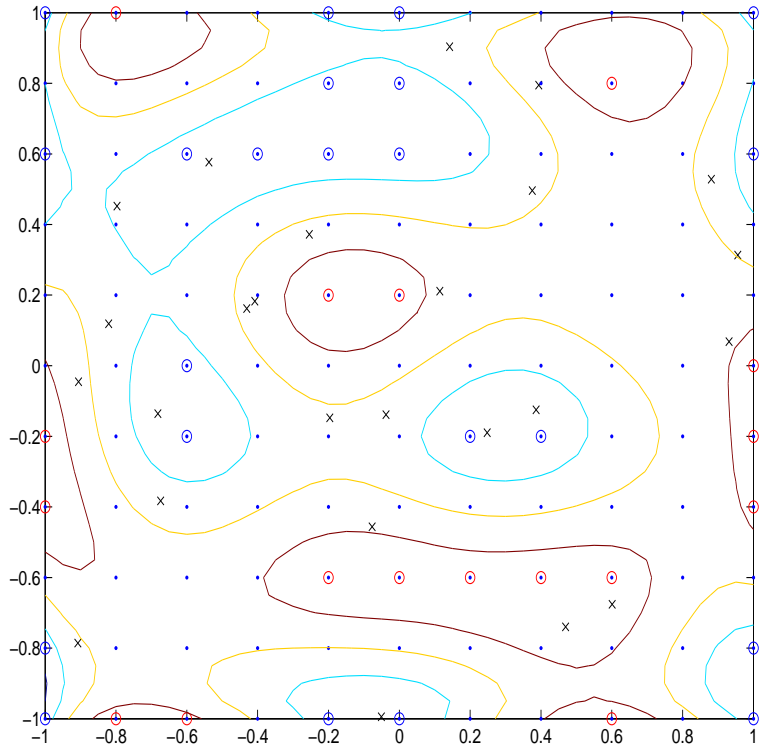}
    \caption{Point sets $X$ (25 crosses), $T$ (121 dots), and extremal points
      (39 circles around points of $T$),
      with contours of the approximation error. The extremal
      points do not form an $H$-set. Signs are indicated by blue or red circles.
      \RSlabel{FigLag}}
\end{center} 
\end{figure}
\biglf
Ignoring what the optimizer says, and aiming at a smaller $\mu$,
one can go for all
points with errors above  $\mu=0.0760$, for instance. This yields only 23
points,
see Figure \RSref{Figmin},
and these do not form an $H$-set either. 
One might argue that
$N=23$ is too small for $n=25$ to make an $H$-set possible, but here and in
other examples on regular points one has
dependent homogeneous equations for the $H$-set
condition \eref{eqdual}, reducing the degrees of freedom.
\begin{figure}[hbt]
  \begin{center}
    \includegraphics[width=10.0cm,height=10.0cm]{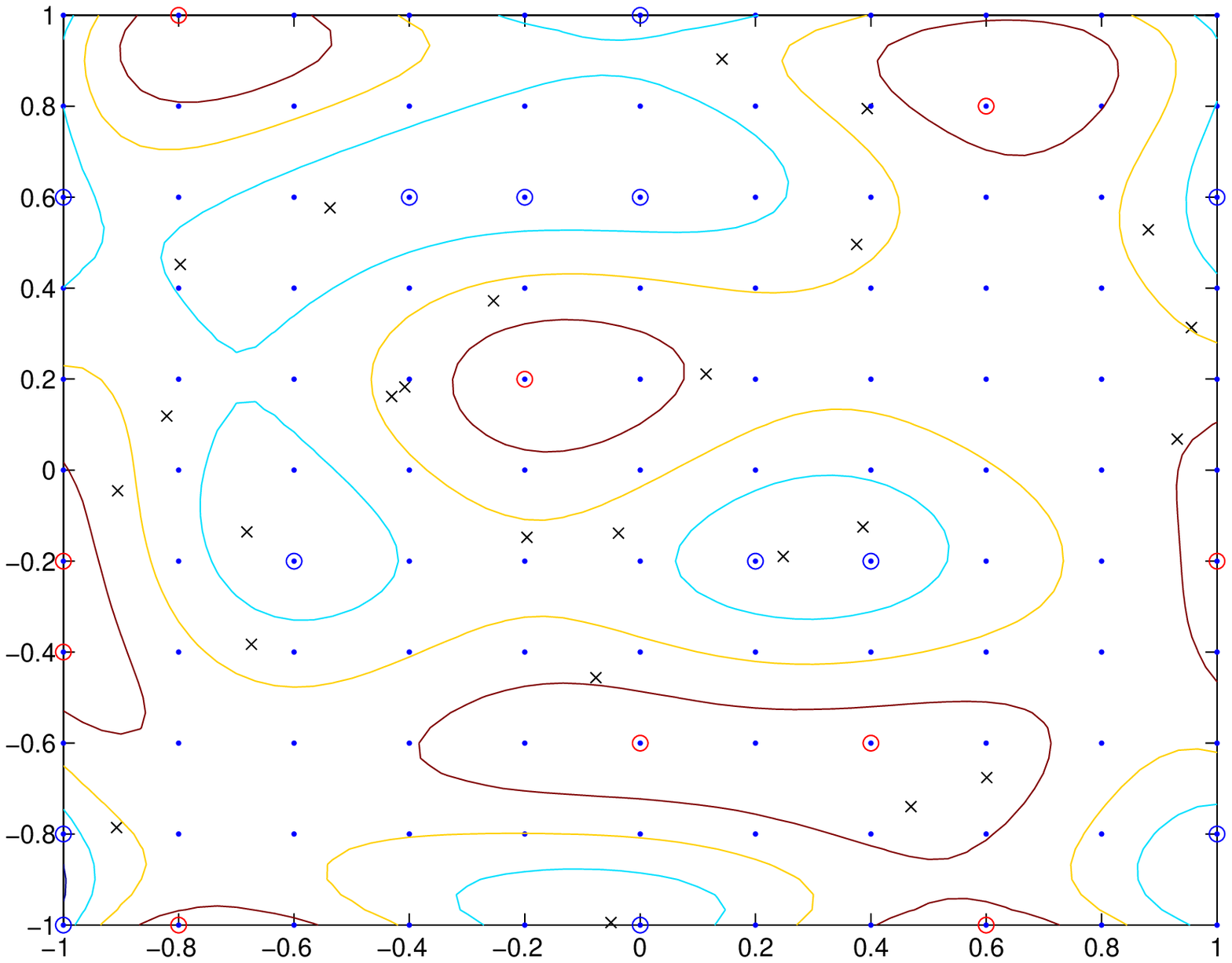}
    \caption{Point sets $X$ (25 crosses), $T$ (121 dots), and extremal points
      (23 circles around points of $T$),
      with contours of the approximation error. Two homogeneous conditions were
      dependent from the others, leading to a 23x23 situation. The extremal
      points do not form an $H$-set.
      \RSlabel{Figmin}}
\end{center} 
\end{figure}
\biglf
But one may take even more points,
by allowing smaller $\mu$ and getting more degrees
of freedom for the $H$-set, by admitting all points that have an absolute error
of $\mu$ or more. It turns out that one has to go down to $\mu=0.0077$
to get an $H$ set of 112 points, see Figure \RSref{Figmax}. But for large $H$,
the maximization of $1^Tw$ shifts large weights to fewer components, 
and thus the set $H$ can be reduced by skipping the zero
components. See Figure \RSref{Figreduced} showing the reduction from 112 to 27
points. Unfortunately, this reduction does not improve $\mu$ reasonably, because
it does not select peak points. It works on coefficients, not on values.
\begin{figure}[hbt]
  \begin{center}
    \includegraphics[width=11.0cm,height=11.0cm]{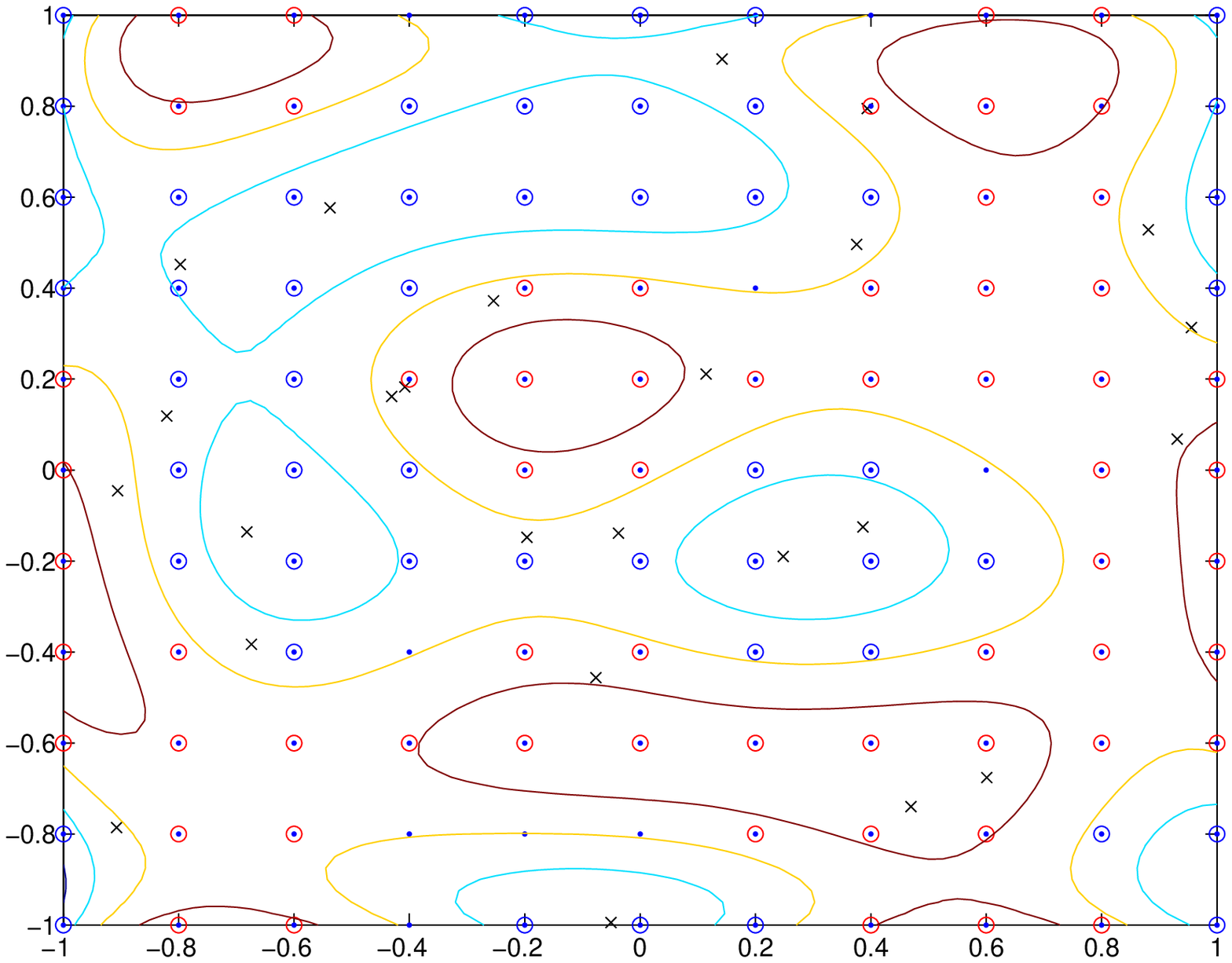}
    \caption{Point sets $X$ (25 crosses), $T$ (121 dots), and $H$-set 
      (112 circles around points of $T$).
      \RSlabel{Figmax}}
\end{center} 
\end{figure}
\begin{figure}[hbt]
  \begin{center}
    \includegraphics[width=11.0cm,height=11.0cm]{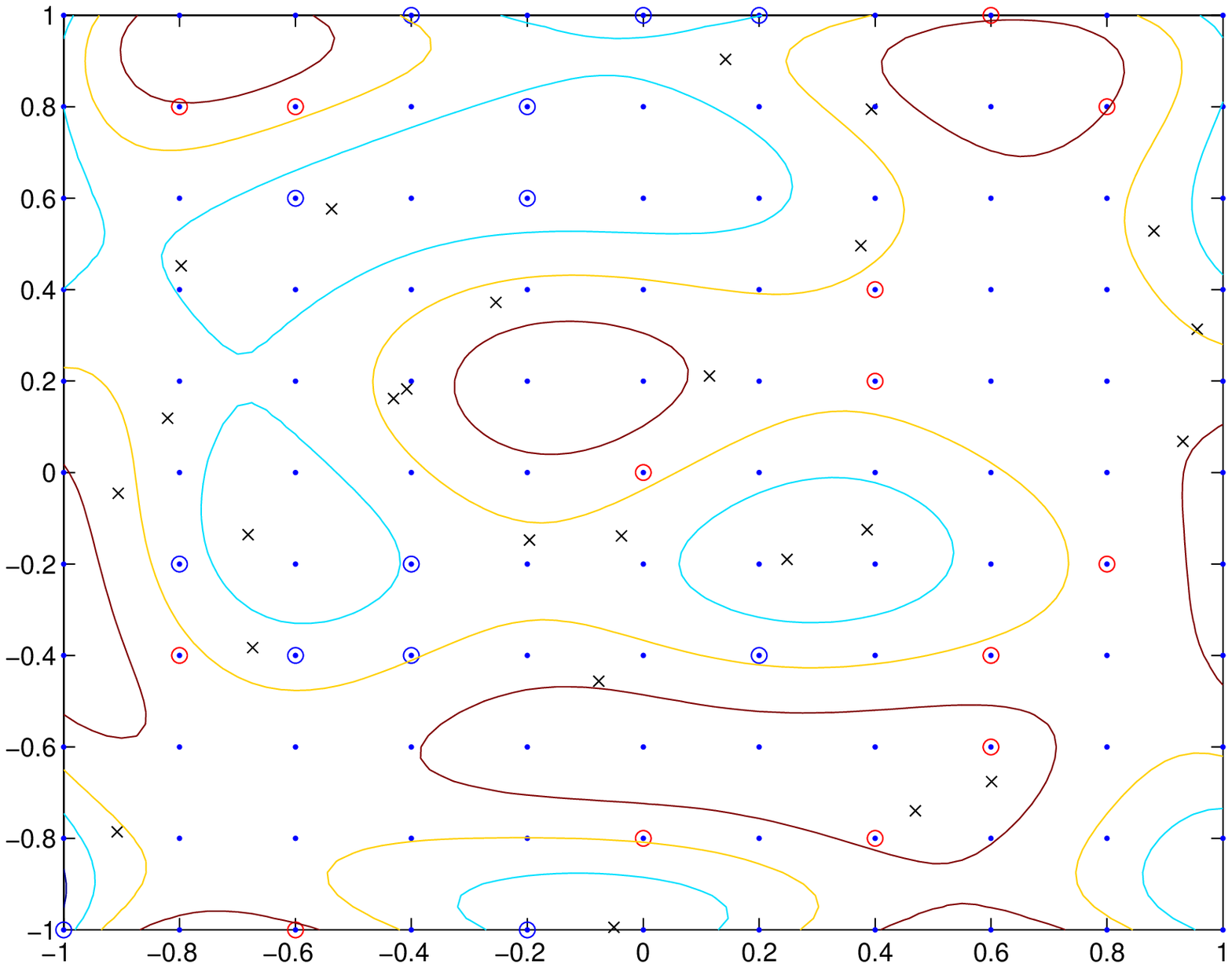}
    \caption{Point sets $X$ (25 crosses), $T$ (121 dots), and reduced
      $H$ set 
      (27 circles around points of $T$)
      \RSlabel{Figreduced}}
\end{center} 
\end{figure}
\section{Kernel-Based Divided Differences}\RSlabel{SecKBDD}
We now consider the case $T=X\cup\{\xi\}$ with $\xi\notin X$
that works perfectly fine
for univariate polynomial approximation, leading to alternation and
divided differences. Generically,  Chebyshev approximation by an $n$-dimensional
space on a set of $n+1$ points should lead to ``equioscillation'', i.e.
the optimal error $\eta^*$ should be
attained at all $n+1$ points, with different
signs. But this  cannot be expected in multivariate situations,
and here we check the case of kernel-based trial spaces.
\biglf
We go into the dual situation and apply 
existence and uniqueness of kernel-based interpolants to get
the unique function $g_\xi\in V_{X\cup\{\xi\}}$
that vanishes on $X$ and is one at $\xi$. If we generally denote
the Lagrangian with respect to a point $y\in Y$ and based on $Y$ as
$u_y^Y$, the function $g_\xi$ is the Lagrangian
$u_\xi^{X\cup\{\xi\}}$ and can be written as 
$$
\dfrac{1}{P_X^2(\xi)} \left(K(x,\xi)-\sum_{i=1}^n u_{x_i}^X(x)K(x_i,\xi)\right)
$$
due to 
$$
  K(\xi,\xi)-\sum_{i=1}^n u_{x_i}^X(\xi)K(x_i,\xi)=P_X^2(\xi)
$$
by definition of the
Power Function $P_X$.
The Lagrangians on $X$ have the form
$$
u_{x_i}^X(x)=\sum_{j=1}^n\alpha^X_{ij}K(x,x_j)
$$
with the $\alpha^X_{i,j}$ being the elements of the
(symmetric) inverse of the kernel matrix based on $X$.
Then the  $\ell_1$ norm of the coefficients of $g_\xi$
in the basis of $V_{X\cup\{\xi\}}$ is obtainable via
$$
\begin{array}{rcl}
  K(x,\xi)-\sum_{i=1}^n u_{x_i}^X(x)K(x_i,\xi)
  &=&
  K(x,\xi)-\sum_{i=1}^n \sum_{j=1}^n \alpha^X_{ij}K(x,x_j)K(x_i,\xi)\\
  &=&
  K(x,\xi)-\sum_{j=1}^n K(x,x_j) \sum_{i=1}^n \alpha^X_{ij}K(x_i,\xi)\\
  &=&
  K(x,\xi)-\sum_{j=1}^n K(x,x_j) \sum_{i=1}^n \alpha^X_{ji}K(x_i,\xi)\\
  &=&
  K(x,\xi)-\sum_{j=1}^n K(x,x_j)u_{x_j}^X(\xi)\\
\end{array} 
$$
as
$$
\dfrac{1}{P_X^2(\xi)} \left(1+ \sum_{j=1}^n |u_{x_j}^X(\xi)| \right)
=\dfrac{1+L_X(\xi)}{P_X^2(\xi)}
$$
using the definition of the Lebesgue function $L_X$.
The solution vector $w^*_T$ for the dual problem thus is unique
and
has coefficients
$$
\begin{array}{rcl}
\dfrac{1}{1+L_X(\xi)} & & \hbox{  at } \xi\\
\dfrac{-u_{x_i}^X(\xi)}{1+L_X(\xi)} & & \hbox{  at } x_i\in X,
\end{array}
$$
up to a fixed sign, because the Power Function cancels out.
Using the standard interpolant $s_{X,f}$ to $f$ on $X$ in its
Lagrange representation, and ignoring a possible sign of $w^*_T$, we find 
$$
  f_T^Tw^*_T =
  \dfrac{f(\xi)-s_{X,f}(\xi)}{1+L_X(\xi)},
$$
  and this is the analog of the divided difference
  in the context of discrete Chebyshev approximation on
  $n+1$ points. In fact, its absolute value
\bql{eqabsval}
\begin{array}{rcl}
  |f_T^Tw^*_T| &=&
  \dfrac{|f(\xi)-s_{X,f}(\xi)|}{1+L_X(\xi)}\\
&=& \eta^*(f,X\cup\{\xi\})
\end{array}
\eq
determines the
maximal error $\eta^*(f,X\cup\{\xi\})$
for the best discrete approximation $s^*_{X,\xi,f}$
to $f$ from the space $V_X$ on $X\cup\{\xi\}$,
because there is no duality gap and $w_T^*$ can only change by its sign.
The complementary slackness
conditions finally produce an $H$-set consisting of $\xi$
and the points $x_i$ of $X$ for which $u_{x_i}(\xi)$ is nonzero.
These must be extremal points, and the sign there is the sign
of $u_{x_i}(\xi)$. Note that \RScite{mueller-schaback:2009-1}
has a similar notion of divided differences in context with Newton bases.
\biglf
In the polynomial case, all
Lagrangians must be nozero at additional points
due to the Fundamental Theorem of Algebra, must change signs
between zeros, 
and therefore one has alternation on all points of $X\cup\{\xi\}$.
In the kernel case, the absolute
errors in all $n+1$ points are equal as long as $\xi$  does not lie
on a zero set of one of the Lagrangians $u_{x_j}^X$. This may be called the
``nondegenerate'' situation of full equioscillation,
if degeneration counts the number
of points where the error is not extremal.
Generically, through each
$x_j\in X$ there will be $n-1$ zero sets defined by the other Lagrangians.
See Figure \RSref{FigExcCur} for the case of the numerical example of Section
\RSref{SecNumExp} using 25 scattered points in $[-1,+1]^2$.
If $\xi$ does not hit one of the curves, there will be no degeneration,
and if $\xi$ moves over the  zero curve of $u_{x_j}^X$, the sign
of the error at $x_j$ will swap.  In view of multiple intersections,
the orders of degeneration
may vary, but with probability
one there is no degeneration, if $\xi$ is sampled uniformly over $[-1,+1]^2$. 

\begin{figure}[hbt]
  \begin{center}
    \includegraphics[width=10.0cm,height=10.0cm]{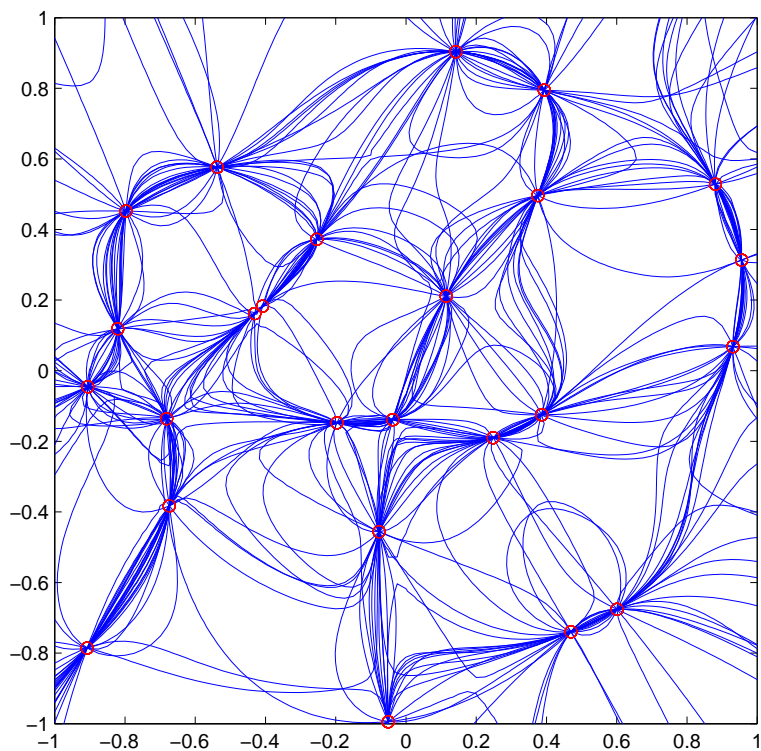}
    \caption{Zero sets of Lagrangians for the 25 points of $X$ (red circles)
      using the Gaussian at scale one.
      \RSlabel{FigExcCur}}
\end{center} 
\end{figure}
\biglf
Figure \RSref{FigDivDiff} shows the divided difference as a function of
$\xi\in[-1,+1]^2$, while Figure \RSref{FigZeroSet} shows the
zero set of the standard interpolation error. 
Note that the points of the zero set can
be added to $X$ without changing the interpolant. This means that the usual
error bounds in terms of fill distances
$$
h(X,\Omega):=\sup_{y\in \Omega}\;\;\min_{x\in X}\|x-y\|_2
$$
should be replaced by the $f$-dependent quantity
$$
\sup_{y\in \Omega}\;\;\inf \{  \|x-y\|_2  \;: \; f(x)=s_{X,f}(x)\}\leq h(X,\Omega).
$$
The $f$-greedy
point selection strategy of \RScite{mueller:2009-1} works similarly,
but picks extrema of the current interpolation error
$f-s_{X,f}$, not points of largest distance to the
zero set. It could as well be changed to pick the
point $\xi$ where the right-hand side of
\eref{eqabsval} is maximal. These variations are open for further research.

\begin{figure}[hbt]
  \begin{center}
    \includegraphics[width=10.0cm,height=10.0cm]{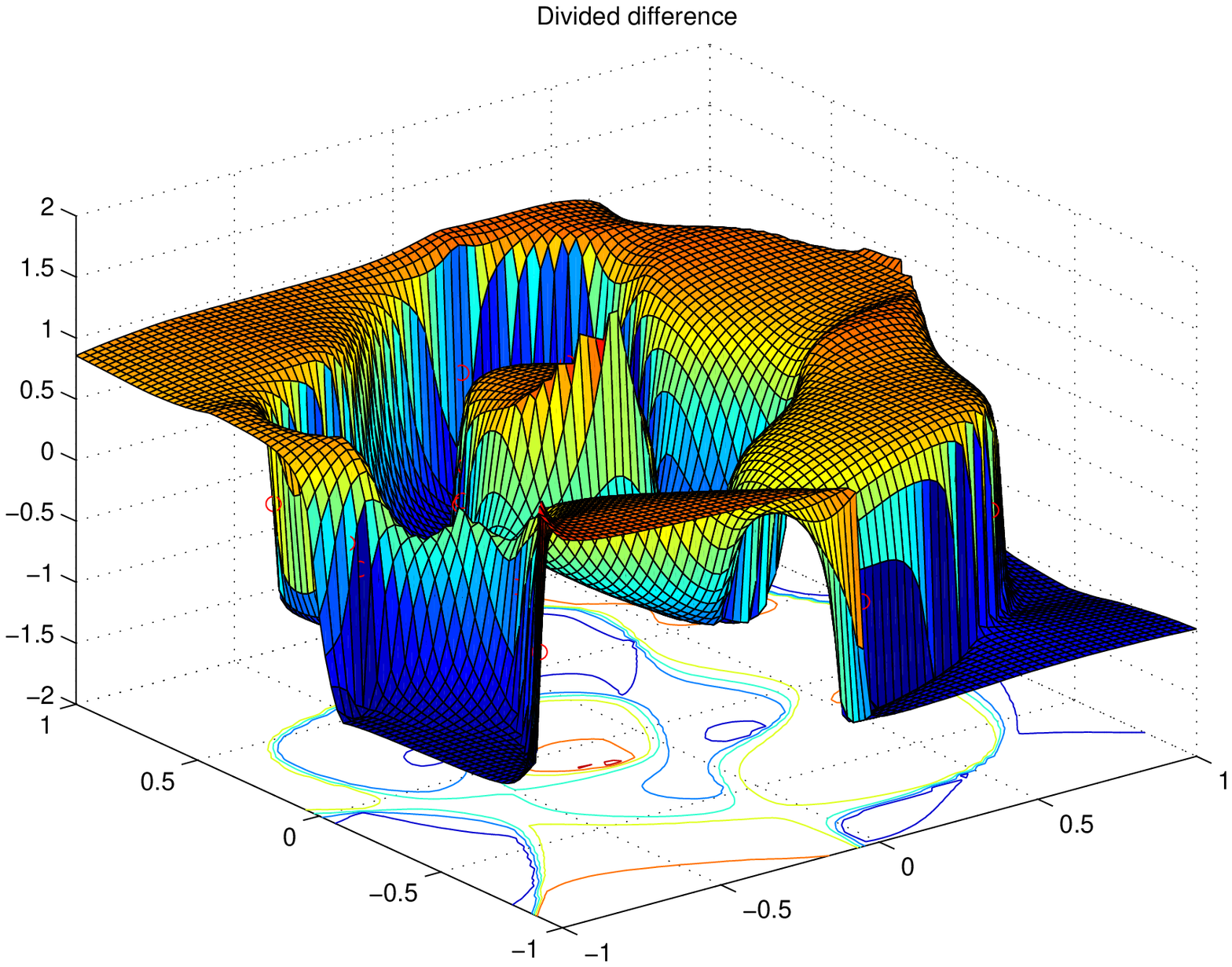}
    \caption{Divided difference as a function of $\xi$
      for the 25 points of $X$ 
      using the Gaussian at scale one and approximating the peaks function.
      \RSlabel{FigDivDiff}}
\end{center} 
\end{figure}

\begin{figure}[hbt]
  \begin{center}
    \includegraphics[width=10.0cm,height=10.0cm]{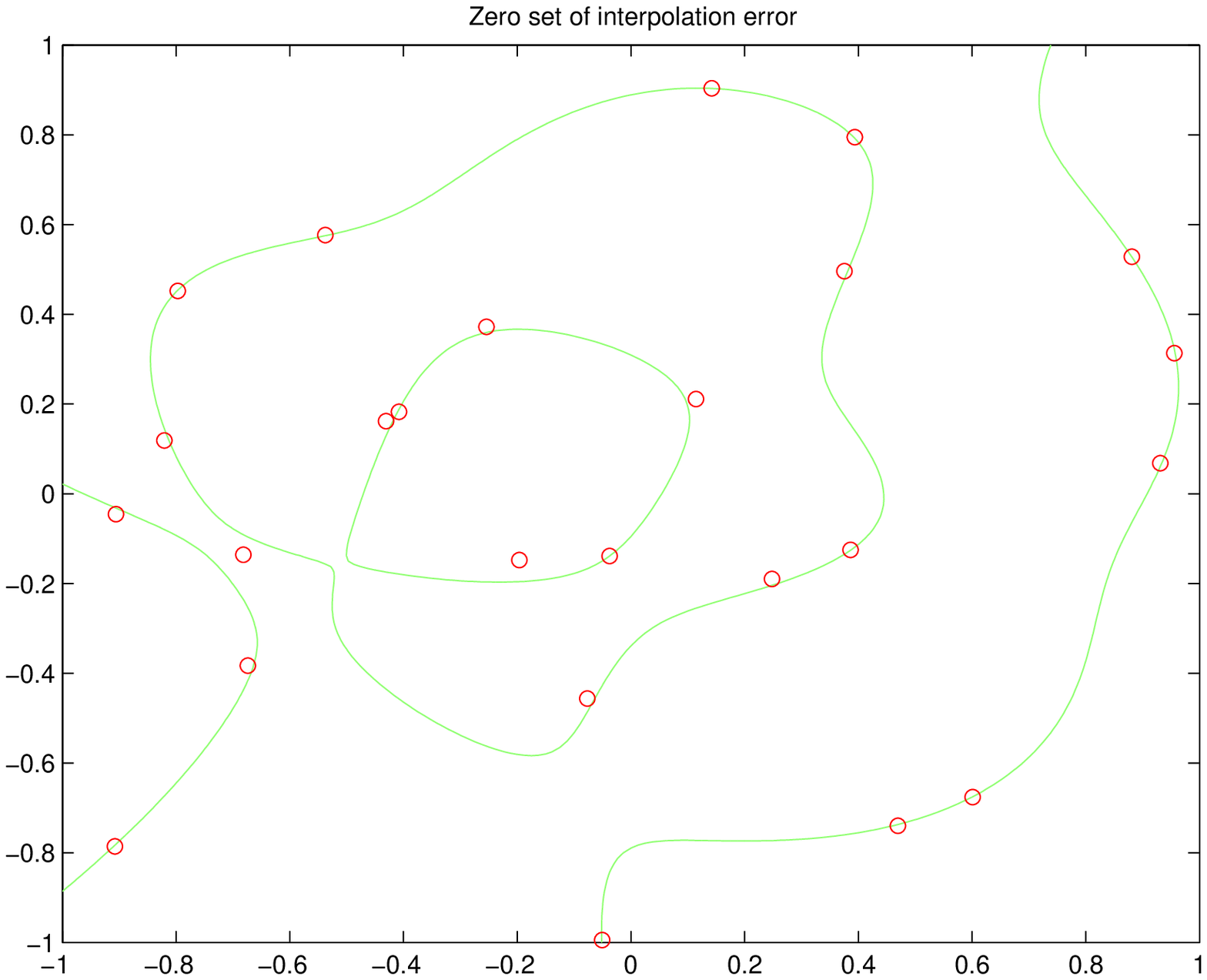}
    \caption{Zero set of the interpolation error
      on the 25 points of $X$ (red circles)
      using the Gaussian at scale one and interpolating the peaks function.
      \RSlabel{FigZeroSet}}
\end{center} 
\end{figure}
\biglf
There is not much known about what happens for
interpolation or approximation using unsymmetric kernel matrices
with entries $K(t_k,x_j),\;1\leq k\leq N,\;1\leq j\leq n$. The above case
$N=n+1$ with $T=X\cup\{\xi\}$ is a first step. 

\bibliographystyle{plain}

\end{document}